\documentclass[a4paper]{article}
\usepackage{amssymb,amsmath}
\newtheorem{theo}{\sc Theorem}[section]
\newtheorem{prop}[theo]{\sc Proposition}
\newtheorem{rema}[theo]{\sc Remark}
\newtheorem{exer}[theo]{\sc Exercise}

\begin{document}

\title{Spherical image of the local Hecke series of genus four}

\author{Kirill Vankov\\
\small Institut Fourier, Universit\'e Grenoble 1\\
\small UFR de Math\'ematiques, UMR 5582\\
\small BP 74, 38402 Saint-Martin d'H\`eres Cedex\\
\small France\\
\small e-mail: kvankov$@$fourier.ujf-grenoble.fr\\
\small phone: +33 476514808
}
\date{}
\maketitle

\begin{abstract}
We obtain an explicit formula for the spherical image of the polynomial fraction conjectured by Shimura in 1963 for generating Hecke series in particular case of genus 4.  As in our previous work we used the Satake spherical map for ${\rm Sp}_4$ and formulas by Andrianov.  Key words: Hecke operators, Dirichlet series, Euler factorization, spherical map.
\end{abstract}

\section{Introduction}

In this article we continue computations of formal power series in Hecke algebra.  The previous article \cite{PaVa06} contains the necessary theory and formulas we used in general case of genus $n$ and the explicit result for $n=3$.

Consider the group of positive symplectic similitudes
\begin{align}
\nonumber
&\mathrm{S}=\mathrm{S}^n={\rm GSp}_n^+(\mathbb{Q})=\{M\in \mathrm{M}_{2n}(\mathbb{Q})\ |\ {}^tMJ_nM=\mu(M)J_n, \mu(M)>0\} \ ,
\\
\nonumber
\mbox{where }& J_n = \left( \begin{array}{cc}\mathbf{0}_n& \mathbf{I}_n \\ -\mathbf{I}_n& \mathbf{0}_n\end{array}\right) . 
\end{align}

For the Siegel modular group $\Gamma={\rm Sp}_n(\mathbb{Z})$ consider the double cosets $(M)=\Gamma M \Gamma \subset \mathrm{S}$ and the Hecke operators $\mathbf{T}(a)= \sum_{M\in \mathrm{SD}_n(a)}(M)$, where $M$ runs through the folowing integer matrices
\begin{align}
\nonumber
&\mathrm{SD}_n(a)=\{{\rm diag} (d_1, \dots, d_n;e_1, \dots, e_n)\ |\ d_i| d_{i+1}, d_n| e_{n}, e_{i+1}|e_i, d_ie_i=a\}. 
\end{align}

Let us use the notation for the Hecke operators
\begin{align}
\nonumber
&\mathbf{T}(d_1, \dots, d_n;e_1, \dots, e_n))=({\rm diag}(d_1, \dots, d_n;e_1, \dots, e_n)).
\end{align}

Let $p$ be a prime.  The formal generating series of Hecke operators of genus $n$ for the symplectic group is written as 
\begin{equation*}
\mathbf{D}_p(\mathsf{X})=\sum_{\delta=0}^{\infty}\mathbf{T}
(p^{\delta})\mathsf{X}^{\delta} \in {\mathcal L}_{n,\mathbb{Z}}
[\![\mathsf{X}]\!] \, ,
\end{equation*}
where the Hecke algebra ${\mathcal L}_{n,\mathbb{Z}}=\mathbb{Z}[\mathbf{T}(p),\, \mathbf{T}_1(p^2),\,\dots\, ,\mathbf{T}_n(p^2)]$ is generated by the following Hecke operators:
\begin{itemize}
\item[] $\mathbf{T}(p)= \mathbf{T}
  (\underbrace{1,\dots,1}_{n},\underbrace{p,\dots,p}_{n})$ \, and
\item[] $\mathbf{T}_i(p^2)= \mathbf{T}
  (\underbrace{1,\dots,1}_{n-i},\underbrace{p,\dots,p}_{i},
  \underbrace{p^2,\dots,p^2}_{n-i},\underbrace{p, \dots, p}_{i})$,
  \mbox{ for } $i=1,\dots,n$.
\end{itemize}

Applying the spherical map $\Omega$ we can carry out all calculations in polynomial rings instead of the Hecke rings of the symplectic group.  This theory is developed in \cite{Sh63}, \cite{An87} and \cite{AnZh95}.  Our previous article \cite{PaVa06} describes in details the method of Andrianov for finding images of Hecke operators.  We adopted the formulas of the article \cite{An70}.  Result is presented in the form of $sym_{i_1,i_2,i_3,i_4}$ polynomials, which are the same symmetrical polynomials as in \cite{PaVa06}, but depend on 4 variables $x_1$, $x_2$, $x_3$ and $x_4$.  Formal algebraic computations were carried out in Maple as for the case $n=3$.  For $n=4$ the number of $\omega(t)$ images increased from 28 ($n=3$) to 680.  It took hours of processor time to obtain the final result.  Intermediate polynomials 
used megabytes of disk space.

\section{The spherical image of $\mathbf{D}_p$}

\begin{theo}
\label{main_result}
 For $n=4$ there is the following explicit polynomial presentation: 
\begin{equation*}
\Omega(\mathbf{D}_p(\mathsf{X})) = \frac{P_4(\mathsf{X})}{Q_4(\mathsf{X})}\;,
\end{equation*}
where
\begin{equation*}
\begin{split}
Q_4(\mathsf{X})=
&(1-x_0\mathsf{X})(1-x_0x_1\mathsf{X})(1-x_0x_2\mathsf{X})(1-x_0x_3\mathsf{X})(1-x_0x_4\mathsf{X})\times \\
&\times (1-x_0x_1x_2\mathsf{X})(1-x_0x_1x_3\mathsf{X})(1-x_0x_1x_4\mathsf{X})(1-x_0x_2x_3\mathsf{X})\times \\
&\times (1-x_0x_2x_4\mathsf{X})(1-x_0x_3x_4\mathsf{X})(1-x_0x_1x_2x_3\mathsf{X})(1-x_0x_1x_2x_4\mathsf{X})\times \\
&\times (1-x_0x_1x_3x_4\mathsf{X})(1-x_0x_2x_3x_4\mathsf{X})(1-x_0x_1x_2x_3x_4\mathsf{X})
\end{split}
\end{equation*}
and
\begin{equation*}
P_4(\mathsf{X})= \sum_{k=0}^{14} K_k(p,x_0,x_1,x_2,x_3,x_4)\mathsf{X}^k
\end{equation*}
with coefficients $K_k$ are listed in the following table.
\end{theo}

\noindent
\begin{tabular}{|p{12cm}|}
\hline 
\\
$K_{0} = 1$ \\
\\
\hline 
\\
$K_{1} = 0$ \\ 
\\ 
\hline 
\end{tabular}
\begin{tabular}{|p{12cm}|}
\hline 
\\
\begin{math}
K_{2} = -{\displaystyle\frac{x_0^2}{p^2}} \times\left(\begin{array}{l}
p\,          (sym_{2211}+sym_{2110}+sym_{1100})+ \\
(p^2+p+1)    (sym_{2111}+sym_{1110})+ \\
(2p^2+4p+1)\, sym_{1111} 
\end{array}
\right)
\end{math}\\ 
\\
\hline 
\\
\begin{math}
K_{3} = {\displaystyle\frac{x_0^3}{p^3}} \times 
\left(
\begin{array}{l}
(p^2+p)(sym_{3222}+sym_{3221}+sym_{3211}+sym_{3111}+sym_{2220}+\\
 \ \ \ \ \ sym_{2210}+sym_{2110}+sym_{1110}) + \\
(p^3 + 5p^2 + 5p + 1)(sym_{2222} + sym_{2221} + sym_{2211} + sym_{2111} + \\
 \ \ \ \ \ sym_{1111})
\end{array}
\right)
\end{math}\\
\\
\hline 
%
\\
\begin{math}
K_{4} = -{\displaystyle\frac{x_0^4}{p^4}} \times
\left(
\begin{array}{l}
p^2 (sym_{4322}+sym_{4221}+sym_{3220}+sym_{2210}) + \\
p\,(p^2+p+1)(sym_{4222}+sym_{3333}+sym_{3331}+sym_{3311}+\\
 \ \ \ \ \ sym_{3111}+sym_{2220}+sym_{1111}) + \\
p\,(p^2+4p+1)(sym_{3332}+sym_{3321}+sym_{3211}+sym_{2111}) + \\
p\,(3p^2+6p+4)(sym_{3322}+sym_{3221}+sym_{2211}) + \\
(5p^3+15p^2+6p+1)(sym_{3222}+sym_{2221}) + \\
(12p^3+22p^2+16p+1)sym_{2222}
\end{array}
\right)
\end{math}\\ 
\\
\hline 
\\
\begin{math}
K_{5} = {\displaystyle\frac{x_0^5}{p^4}} \times
\left(
\begin{array}{l}
(p^2+p)(sym_{4433}+sym_{4432}+sym_{4422}+sym_{4331}+sym_{4321} + \\
 \ \ \ \ \ sym_{4221}+sym_{3311}+sym_{3211}+sym_{2211}) + \\
(4p^2+5p+1)(sym_{4333}+sym_{4332}+sym_{4322}+sym_{4222} + \\
 \ \ \ \ \ sym_{3331}+sym_{3321}+sym_{3221}+sym_{2221}) + \\
(-p^4+14p^2+18p+5)(sym_{3333}+sym_{3332} + \\
 \ \ \ \ \ sym_{3322}+sym_{3222}+sym_{2222})
\end{array}
\right)
\end{math}\\ 
\\
\hline 
\\
\begin{math}
K_{6} = {\displaystyle\frac{x_0^6}{p^6}} \times
\left(
\begin{array}{l}
p^2(p^3-5p-4)(sym_{4432}+sym_{4322}+sym_{3222}+sym_{4443}) + \\
p\,(p^5+5p^4-17p^2-15p-1)(sym_{4333}+sym_{3332}) - \\
p^2(p+1)(sym_{4331}+sym_{3321}+sym_{5332}+sym_{5433}) + \\
p\,(3p^4-12p^2-6p-1)(sym_{4332}+sym_{3322}+sym_{4433}) + \\
p^2(p^3-3p-1)(sym_{4222}+sym_{3331}+sym_{2222}+sym_{4422} + \\
 \ \ \ \ \ sym_{4442}+sym_{4444}+sym_{5333}) + \\
p^3(sym_{6333}-sym_{5443}-sym_{5432}-sym_{5322}-sym_{4431} - \\
 \ \ \ \ \ sym_{4321}+sym_{3330}-sym_{3221}) + \\
(2p^6+12p^5-32p^3-22p^2-4p+1)sym_{3333}
\end{array}
\right)
\end{math}\\ 
\\
\hline 
\\
\begin{math}
K_{7} = -{\displaystyle\frac{x_0^7}{p^5}} \times
\left(
\begin{array}{l}
p\,(p^2-1)(sym_{5544}+sym_{5543}+sym_{5533}+sym_{5442} + \\
 \ \ \ \ \ sym_{5432}+sym_{5332}+sym_{4422}+sym_{4322}+sym_{3322}) + \\
(p^4+4p^3-4p-1)(sym_{5444}+sym_{5443}+sym_{5433} + \\
 \ \ \ \ \ sym_{5333}+sym_{4442}+sym_{4432}+sym_{4332}+sym_{3332}) + \\
(5p^4+14p^3-14p-5)(sym_{4444}+sym_{4443}+sym_{4433} + \\
 \ \ \ \ \ sym_{4333}+sym_{3333})
\end{array}
\right)
\end{math}\\ 
\\
\hline 
\end{tabular}
\begin{tabular}{|p{12cm}|}
\hline 
\\
\begin{math}
K_{8} = {\displaystyle\frac{x_0^8}{p^6}} \times
\left(
\begin{array}{l}
(p^5 + 15p^4 + 17p^3 - 5p - 1)(sym_{5444} +  sym_{4443}) + \\
p^3(-sym_{7444} + sym_{6554} + sym_{6543} + sym_{6433} + sym_{5542} + \\
 \ \ \ \ \ sym_{5432} - sym_{4441} + sym_{4332}) + \\
p\,(4p^3 + 5p^2 - 1)(sym_{5554} + sym_{5543} + sym_{5433} + sym_{4333}) + \\
p^3(p + 1)(sym_{6544} + sym_{6443} + sym_{5442} + sym_{4432}) + \\
p\,(p^3 + 3p^2 - 1)(sym_{4442} + sym_{5333} + sym_{5533} + sym_{5555} + \\
 \ \ \ \ \ sym_{6444} + sym_{5553} + sym_{3333}) - \\
(p^6 - 4p^5 - 22p^4 - 32p^3 + 12p + 2) sym_{4444} + \\
p\,(p^4 + 6p^3 + 12p^2 - 3)(sym_{5544} + sym_{5443} + sym_{4433})
\end{array}
\right)
\end{math}\\ 
\\
\hline 
%
\\
\begin{math}
K_{9} = -{\displaystyle\frac{x_0^9}{p^6}} \times
\left(
\begin{array}{l}
p^2(p+1)(sym_{6655}+sym_{6654}+sym_{6644}+sym_{6553} + \\
 \ \ \ \ \ sym_{6543}+sym_{6443}+sym_{5533}+sym_{5433}+sym_{4433}) + \\
p^2(p^2+5p+4)(sym_{6555}+sym_{6554}+sym_{6544}+sym_{6444} + \\
 \ \ \ \ \ sym_{5553}+sym_{5543}+sym_{5443}+sym_{4443}) + \\
(5p^4+18p^3+14p^2-1)(sym_{5555}+sym_{5554} +  \\
 \ \ \ \ \ sym_{5544}+sym_{5444}+sym_{4444})\\
\end{array}
\right)
\end{math}\\ 
\\
\hline 
\\
\begin{math}
K_{10} = {\displaystyle\frac{x_0^{10}}{p^5}} \times
\left(
\begin{array}{l}
(p^2+p+1)(sym_{4444}+sym_{5553}+sym_{6444}+sym_{6644} + \\
 \ \ \ \ \ sym_{6664}+sym_{6666}+sym_{7555}) + \\
(p^2+4p+1)(sym_{5444}+sym_{6544}+sym_{6654}+sym_{6665}) + \\
p\,(sym_{5543}+sym_{6553}+sym_{7554}+sym_{7655}) + \\
(4p^2+6p+3)(sym_{5544}+sym_{6554}+sym_{6655}) + \\
(p^3+6p^2+15p+5)(sym_{5554}+sym_{6555}) + \\
(p^3+16p^2+22p+12)sym_{5555}\\
\end{array}
\right)
\end{math}\\ 
\\
\hline 
\\
\begin{math}
K_{11} = -{\displaystyle\frac{x_0^{11}}{p^6}} \times
\left(
\begin{array}{l}
(p^2+p)(sym_{7666}+sym_{7665}+sym_{7655}+sym_{7555} + \\
 \ \ \ \ \ sym_{6664}+sym_{6654}+sym_{6554}+sym_{5554}) + \\
(p^3+5p^2+5p+1)(sym_{6666}+sym_{6665}+sym_{6655} + \\
 \ \ \ \ \ sym_{6555}+sym_{5555})\\
\end{array}
\right)
\end{math}\\ 
\\
\hline 
\\
\begin{math}
K_{12} = {\displaystyle\frac{x_0^{12}}{p^6}} \times
\left(
\begin{array}{l}
(p^2+p+1)(sym_{6665}+sym_{7666}) + \\
(p^2+4p+2)sym_{6666} + \\
p\,(sym_{7665}+sym_{7766}+sym_{6655})\\
\end{array}
\right)
\end{math}\\ 
\\
\hline 
\\
$K_{13} = 0$\\
\\
\hline 
\\
\begin{math}
K_{14} = -{\displaystyle\frac{x_0^{14}}{p^6}} \times sym_{7777}
\end{math}\\ 
\\
\hline 
\end{tabular}

\section{Remarks}

We noticed a very interesting symmetry property within the coefficients $K_k$.  Knowing this relation in advance would let to limit computation of coefficients just up to degree 7, skipping the most time consuming higher degree coefficients.

\begin{prop}
Polynomial $P_4(\mathsf{X})$ has the following functional relation between its coefficients $K_k$, $k=0,\dots ,14$:
\begin{equation*}
\begin{split}
 K_{14-k}&(p,x_0,x_1,x_2,x_3,x_4) = \\
 &-p^{-6}(x_0^2x_1x_2x_3x_4)^{7-k}\,
 K_{k}\left(\frac{1}{p},x_0x_1x_2x_3x_4,\frac{1}{x_1},\frac{1}{x_2},\frac{1}{x_3},\frac{1}{x_4}\right)
 \end{split}
\end{equation*}
\end{prop}

\begin{rema}
It is suggested that this functional relation is true for all $n$ in the following form:
\begin{equation*}
 P(x_0, \dots ,x_n,\mathsf{X}) = 
 (-1)^{n-1} \frac{(x_0^2x_1\dots x_n\mathsf{X}^2)^{2^{n-1}-1}}{p^{\frac{n(n-1)}{2}}} P\left( \frac{1}{x_0}, \dots , \frac{1}{x_n}, \frac{p}{\mathsf{X}}\right)
\end{equation*}
\end{rema}

\begin{exer}
The result of the Theorem \ref{main_result} is in a full agreement with the result of the earlier work \cite{PaVa06} for $n=3$ by applying a projection $x_4=0$ (corresponding to Siegel operator acting from ${\rm Sp}_4$ to ${\rm Sp}_3$).
\end{exer}

\bigbreak
The author is very grateful to his academic advisor professor A. A. Panchishkin for posing the problem and active discussions.

\bibliographystyle{amsplain}

\end{document}